# On the Improved Nonlinear Tracking Differentiator based Nonlinear PID Controller Design

Ibraheem Kasim Ibraheem
Electrical Engineering Department
College of Engineering, Baghdad University
Baghdad, Iraq

Wameedh Riyadh Abdul-Adheem
Electrical Engineering Department
College of Engineering, Baghdad University
Baghdad, Iraq

*Abstract*—This paper presents a new improved nonlinear tracking differentiator (INTD) with hyperbolic tangent function in the state-space system. The stability and convergence of the INTD are thoroughly investigated and proved. Through the error analysis, the proposed INTD can extract differentiation of any piecewise smooth nonlinear signal to reach a high accuracy. The improved tracking differentiator (INTD) has the required filtering features and can cope with the nonlinearities caused by the noise. Through simulations, the INTD is implemented as a signal's derivative generator for the closed-loop feedback control system with a nonlinear PID controller for the nonlinear Mass-Spring-Damper system and showed that it could achieve the signal tracking and differentiation faster with a minimum mean square error.

## I. INTRODUCTION

Differentiation of signals in real time is an old and well-known problem. An ideal differentiator would have to differentiate measurement noise with possibly large derivatives along with the signal [1]. In various case studies, the building of a differentiator is inescapable. However, the perfect differentiator could not be synthesized. Without a doubt, together with the principal function, it could differentiate any minor high-frequency noise which is inherent in the signal and may have large derivative values [2].

Designing a differentiator as a single entity is a common design objective for the field of signal processing. The initial procedure is to let some linear dynamical system model to represent the transfer function of the perfect differentiator. Accordingly, the obtained differentiator does not compute the precise derivatives of only noise free signals including the situations when the frequency bandwidth of the signal is limited [2]. Tracking differentiator (TD) design has drawn much consideration in last twenty years because of trailing the high performance of control and navigation system [3].

The traditional high-gain differentiator announced by [4] could follow certain derivatives when the gains lean towards infinity which couldn't feasibly realize. In [2], a sliding mode technique has been used to design differentiator. An upper bound for Lipschitz constant is needed in this kind of differentiator. Nevertheless, the derivative estimation is not soft due to the existence of a discontinuous function. Therefore, a chattering phenomenon occurs in the derivative evaluation. In [5], a universal vigorous precise differentiator has been developed by integrating a sliding modes differentiation with the high-gain differentiator by means of a switching function. Wang Xinhua [6, 7] has suggested a continuous hybridized nonlinear differentiator in which a smattering phenomenon has been decreased adequately. The differentiators in [8] could regularly approach to the correct solution with finite-time exact convergence and start differentiator error [3]. While linear techniques for tracking differentiators design have been adopted in [9, 10].

In [11], two particular high-gain tracking differentiators were proposed. This differentiator was based on the Taylor expansion, the time lagging phenomenon of the traditional high-gain differentiator is reduced effectively.

Also, a fractional order tracking differentiators have been studied recently. In [12], the tracking differentiator was redesigned with fractional-order to provide a fundament for the design of fractional-order ADRC. An Adaptive controller was developed by Wei [13] using fractional order tracking differentiator. In [14], a discrete analog of a fractional order differentiator over Paley–Wiener space are constructed.

As an improvement, a new meta-heuristic optimization algorithm, called cuckoo search algorithm (CSA) was applied by Kumar [15] to determine the optimal coefficients of the finite impulse response-fractional order differentiator

In practice, to achieve high-performance control, many applications based on tracking differentiators have been proposed, such as, the pitch and depth control problem of autonomous underwater vehicle (AUV) in diving plane [16], detection of harmonic current in single phase active power filters [17], geomagnetic attitude detection systems [18], the position and speed detection system as well as suspension system of maglev train [19], electric vehicles (EVs) [20], etc.

In this work, an improved tracking differentiator is proposed, and its stability is tested based on Lyapunov technique. The peaking phenomenon is presented through time domain analysis, while frequency domain analysis proves that the proposed nonlinear tracking differentiator attenuates signals with a certain frequency band.

The paper is organized as follows: in Section II an improved nonlinear tracking differentiator is proposed, and the main convergence results are presented. Section III explains using the INTD with nonlinear PID controller. The mathematical model of the nonlinear Mass-Spring-Damper is introduced in Section IV. The numerical results are presented in Section V to verify the effectiveness of the proposed INTD. Finally, the conclusions are provided in section VI.

## II. THE IMPROVED NONLINEAR TRACKING DIFFERENTIATOR (INTD)

The enhanced nonlinear tracking differentiator based on the hyperbolic tangent function is proposed as follows:

$$\left.\begin{array}{l}\dot{z}_1 = z_2 \\ \dot{z}_2 = -R^2 \tanh\left(\frac{\beta z_1 - (1-\alpha)v}{\gamma}\right) - Rz_2\end{array}\right\} \quad (1)$$

Where $z_1$ is tracking the input $v$, and $z_2$ tracking the derivative of input $v$. the parameters $\alpha, \beta, \gamma, and\ R$ are appropriate design parameters, where $0 < \alpha < 1, \beta > 0, \gamma > 0, and\ R > 0$.

**Lemma 1: (Convergence of the INTD system):** the improved tracking differentiator described by (1) with its design parameters is globally asymptotically stable.

**Proof:** Let us assign $V_l(\mathbf{z}) = R\frac{\gamma}{\beta}\ln\cosh\left(\frac{\beta z_1}{\gamma}\right) + \frac{1}{2}z_2^2$ as a Lyapunov function to system (1). Where $V_l(\mathbf{z}) > 0$ if and only if $\mathbf{z} \neq 0$, and $V_l(\mathbf{z}) = 0$ if and only if $\mathbf{z} = 0$

Now,

$$\dot{V}_l(\mathbf{z}) = -Rz_2^2 \text{ and}$$

$$\dot{V}_l(\mathbf{z}) \leq 0 \text{ for all } z_2$$

This leads to $\dot{V}_l(\mathbf{0}) = 0$ at the origin by Lasalle's theorem[21]. Since $V_l(\mathbf{z}) \to \infty$ for $\|z\| \to \infty$, then the system is globally asymptotically stable (GAS). □

**Lemma 2: (Arrival phase):** consider the system (1) if $\frac{\beta z_1 - (1-\alpha)v}{\gamma} \gg 1$; then $\forall t > 0$, the term $\frac{\beta z_1 - (1-\alpha)v}{\gamma}$ will be decreased until it reaches the tracking phase where $\left|\frac{\beta z_1 - (1-\alpha)v}{\gamma}\right| \ll 1$.

**Proof:** Since $\frac{\beta z_1 - (1-\alpha)v}{\gamma} \gg 1$, Then $\tanh\left(\frac{\beta z_1 - (1-\alpha)v}{\gamma}\right) \to 1$

So that

$$\left.\begin{array}{l}\dot{z}_1 = z_2 \\ \dot{z}_2 = -Rz_2 - R^2\end{array}\right\} \quad (2)$$

The solution of system (2) with the initial condition $\mathbf{z}(0) = [z_1(0)\ z_2(0)]^T$ is given as

$$z_1(t) = -Rt - \left(1 + \frac{z_2(0)}{R}\right)e^{-Rt} + z_1(0) + \frac{z_2(0)}{R} + 1$$

$$z_2(t) = -R + (R + z_2(0))e^{-Rt}$$

If $z_2(0) = 0$, then $z_1(t)$ is a decreasing function for $t > 0$ until it reaches the tracking phase where $\left|\frac{\beta z_1 - (1-\alpha)v}{\gamma}\right| \ll 1$. □

**Corollary 1:** for the system given by (1) if $\frac{\beta z_1(t) - (1-\alpha)v(t)}{\gamma} \ll -1$, then $\forall t > 0$, the term $\frac{\beta z_1(t) - (1-\alpha)v(t)}{\gamma}$ will be increased until the system reaches the tracking phase where $\left|\frac{\beta z_1 - (1-\alpha)v}{\gamma}\right| \ll 1$.

**Proof:** By the same way of lemma (2), for $\frac{\beta z_1(t) - (1-\alpha)v(t)}{\gamma} \ll 1$, and $z_2(0) = 0$, then $z_1(t)$ increasing for $t > 0$ until it reaches the tracking phase where $\left|\frac{\beta z_1 - (1-\alpha)v}{\gamma}\right| \ll 1$. □

**Lemma 3: (tracking phase):** consider system (1) for $\left|\frac{\beta z_1 - (1-\alpha)v}{\gamma}\right| \ll 1$ then both tracking error $e_t(t) = v(t) - \frac{\beta}{1-\alpha}z_1(t)$, and differentiation error $e_d(t) = \dot{v}(t) - \frac{\beta}{1-\alpha}z_2(t)$ tends to zero for finite input signal.

**Proof:** Since $\frac{\beta z_1 - (1-\alpha)v}{\gamma} \ll 1$, Then $\tanh\left(\frac{\beta z_1 - (1-\alpha)v}{\gamma}\right) \to \left(\frac{\beta z_1 - (1-\alpha)v}{\gamma}\right)$. So that,

$$\left.\begin{array}{l}\dot{z}_1 = z_2 \\ \dot{z}_2 = -R^2\left(\frac{\beta z_1 - (1-\alpha)v}{\gamma}\right) - Rz_2\end{array}\right\} \quad (3)$$

Taking Laplace transform to (3), then

$$\begin{bmatrix}Z_1(s)\\Z_2(s)\end{bmatrix} = \begin{bmatrix}\frac{\frac{R^2(1-\alpha)}{\gamma}}{s^2+Rs+\frac{R^2\beta}{\gamma}}\\\frac{\frac{R^2(1-\alpha)s}{\gamma}}{s^2+Rs+\frac{R^2\beta}{\gamma}}\end{bmatrix}V(s) \quad (4)$$

The tracking error associated with the tracking phase is

$$e_t(t) = v(t) - \frac{\beta}{1-\alpha}z_1(t)$$

$$E_t(s) = V(s) - \frac{\beta}{1-\alpha}Z_1(s)$$

The transfer function of the tracking error w.r.t input $v$ is given as

$$L_t(s) = \frac{E_t(s)}{V(s)} = \frac{s(s+R)}{s^2+Rs+\frac{R^2\beta}{\gamma}}$$

So that,

$$l_t(\infty) = \lim_{s\to 0}sL_t(s) = 0 \quad (5)$$

The differentiation error associated during the tracking phase is described as

$$e_d(t) = \dot{v}(t) - \frac{\beta}{1-\alpha}z_2(t)$$

$$E_d(s) = sV(s) - \frac{\beta}{1-\alpha}Z_2(s)$$

The transfer function of the differentiation error w.r.t the input derivative is

$$L_d(s) = \frac{E_d(s)}{sV(s)} = \frac{s(s+R)}{s^2+Rs+\frac{R^2\beta}{\gamma}}$$

So that,

$$l_d(\infty) = \lim_{s \to 0} sL_d(s) = 0 \quad (6)$$

Therefore (5) and (6) complete the proof. □

**Theorem 1:** Consider system (1), then for any value of $\left|\frac{\beta z_1 - (1-\alpha)v}{\gamma}\right|$,

$$\lim_{t \to \infty} \left|\frac{\beta z_1(t) - (1-\alpha)v(t)}{\gamma}\right| = 0$$

and

$$\lim_{t \to \infty} \left|\frac{\beta z_2(t) - (1-\alpha)\dot{v}(t)}{\gamma}\right| = 0$$

**Proof:**

By using Lemma (2) and (3). □

**Lemma 4: (Time domain analysis):** Consider the system (1) which satisfies (4). If $\beta \gg 1, 0 < \gamma < 1, R \gg 1$, with $0 < \alpha < 1$, then the system (1) has a high undamped natural frequency, a small damping ratio, and peaking phenomenon.

**Proof:**

$$\frac{Z_2(s)}{sV(s)} = \frac{\frac{R^2(1-\alpha)}{\gamma}}{s^2 + Rs + \frac{R^2\beta}{\gamma}} = \left(\frac{1-\alpha}{\beta}\right)\frac{\omega_n^2}{s^2 + 2\xi\omega_n s + \omega_n^2}$$

Where,

$\omega_n = R\sqrt{\frac{\beta}{\gamma}}$ is the undamped natural frequency (rad/sec)

$\xi = \frac{1}{2}\sqrt{\frac{\gamma}{\beta}}$ is the damping ratio

It's clear that from the values of the parameters $\beta$, $\gamma$, and $R$ that the damping ratio $\xi \ll 1$ implies that the system has an under damped effect which leads to peaking phenomenon. □

**Lemma 5: (The frequency-domain analysis):** Consider a system (1) which satisfies equation (4) with parameters β, γ, and R defined in Lemma 4. The system represents a band-limited differentiator with bandwidth $\omega_n$.

**Proof:** for

$$\frac{Z_2(j\omega)}{V(j\omega)} = \left(\frac{1-\alpha}{\beta}\right)\frac{\omega_n^2 j\omega}{(j\omega)^2 + 2\xi\omega_n j\omega + \omega_n^2}$$
$$= \left(\frac{1-\alpha}{\beta}\right)\frac{j\omega}{\left(\frac{j\omega}{\omega_n}\right)^2 + 2\xi\frac{j\omega}{\omega_n} + 1}$$

if the magnitude of the transfer function $\frac{Z_2(j\omega)}{V(j\omega)}$ is taken as

$$20\log\left|\frac{Z_2(j\omega)}{V(j\omega)}\right| = 20\log\left(\frac{1-\alpha}{\beta}\right) + 20\log\omega - 20\log\sqrt{(1-\left(\frac{\omega}{\omega_n}\right)^2)^2 + \left(2\xi\frac{\omega}{\omega_n}\right)^2}$$

For $\omega \ll \omega_n$ this implies

$$20\log\left|\frac{Z_2(j\omega)}{V(j\omega)}\right| = 20\log\left(\frac{1-\alpha}{\beta}\right) + 20\log\omega$$

Such that $20\log\left(\frac{1-\alpha}{\beta}\right)$ is the correction gain and $20\log\omega$ is the differentiator effect.

On the other hand if $\omega \gg \omega_n$, then

$$20\log\left|\frac{Z_2(j\omega)}{V(j\omega)}\right| = 20\log\left(\frac{1-\alpha}{\beta}\right) + 20\log\omega - 40\log\frac{\omega}{\omega_n}$$

The third term represents the attenuation effect. Therefore, the system has the attenuation effect for $\omega \gg \omega_n$. □

## III. TRACKING DIFFERENTIATOR BASED NONLINEAR PID CONTROLLER

Using nonlinear tracking differentiator, a standard PID controller is transformed into nonlinear PID (NLPID) [22] as shown in Fig. 1. The first tracking differentiator (TD(I)) is used as transient process profile generator, while the second tracking differentiator (TD(II)) is used as state observer to get tracking output $z_1$ and its differential $z_2$. The error, integral, and differential signals are produced by comparing transient process profile to the output of TD(II).

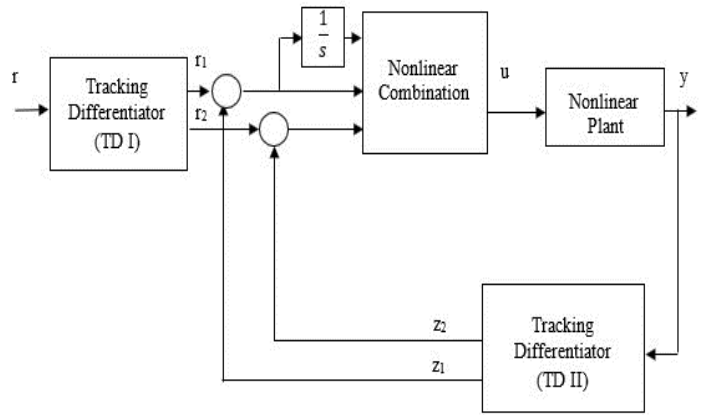

Fig. 1. The traditional structure of NLPID controller

Jing Han has made some investigations on traditional structures and essential properties of nonlinear tracking differentiator. A kind of second-order nonlinear tracking differentiator based on second order bang-bang switch system has been proposed [22]:

$$\dot{x}_1 = x_2$$
$$\dot{x}_1 = -R\, sign\left(x_1 - v(t) + \frac{x_2|x_2|}{2R}\right)$$

where $x_1$ is the desired trajectory and $x_2$ is its derivative. Note that, the parameter $R$ is an application dependent and it is set accordingly to speed up or slow down the transient profile. Then, $x_2$ is denoted as the "tracking differentiator" of $v(t)$.

In order to avoid chattering near the origin, changing the **sign** function to linear saturation function **sat**, then the modified Han TD is represented by:

$$\dot{x}_1 = x_2$$

$$\dot{x}_1 = -R\, sat(x_1 - v(t) + \frac{x_2\,|x_2|}{2R}, \delta)$$

Where

$$sat(A,\delta) = \begin{cases} sign(A), & |A| > \delta \\ \frac{A}{\delta} & |A| \leq \delta \end{cases}$$

The NLPID takes "nonlinear combination" on the three signals. Han [23] proposed the following nonlinear function:

$$fal(e,\alpha,\delta) = \begin{cases} \frac{e}{\delta^{1-\alpha}} & |x| \leq \delta \\ |e|^\alpha sign(e) & |x| \geq \delta \end{cases}$$

The control rule takes:

$$u = \beta_1 fal(e,\alpha_1,\delta_1) + \beta_2 fal(\dot{e},\alpha_2,\delta_2) + \beta_3 fal\left(\int e, \alpha_3, \delta_3\right)$$

Where $\alpha_1, \alpha_2$ and $\alpha_3 \in [0.5\ 1]$

## IV. MATHEMATICAL MODELING OF THE NONLINEAR MASS-SPRING-DUMPER (NMSD) PLANT

A simple nonlinear mass-spring-damper mechanical system as shown in Fig. 2. It is assumed that the stiffness coefficient of the spring, the damping coefficient of the damper and the input term have nonlinearity or uncertainty [23]:

$$M\ddot{x} + g(x,\dot{x}) + f(x) = \varphi(\dot{x})u \quad (7)$$

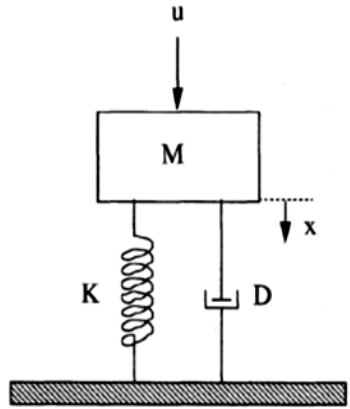

Fig. 2. The nonlinear mass spring damper model

where $M$ is the mass and $U$ is the force, $f(x)$ is the nonlinear or uncertain term with respect to the spring, $g(x,\dot{x})$ is the nonlinear or uncertain term with respect to the damper, and $\varphi(\dot{x})$ is the nonlinear term with respect to the input term.

Assume that $g(x,\dot{x}) = D(c_1 x + c_2 \dot{x}^3), f(x) = c_3 x + c_4 x^3$, and $\varphi(\dot{x}) = 1 + c_5 \dot{x}^3$, assume that $x \in [-a\ a]$, $\dot{x} \in [-b\ b]$, $a, b > 0$. The above parameters are set as follows:

$M = 1.0$, $D = 1.0$, $c_1 = 0.01$, $c_2 = 0.1$, $c_3 = 0.01$, $c_4 = 0.67$, $c_5 = 0$, $a = 1.5$, b = 1.5. Then (6) can be written as:

Then, (7) can be rewritten as follows:

$$\ddot{x} = -0.1\dot{x}^3 - 0.02x - 0.67x^3 + u \quad (8)$$

The state space representation of the nonlinear mass-spring-dumper model is:

$$\left.\begin{array}{l} \dot{x}_1 = x_2 \\ \dot{x}_2 = -0.1 x_2^{\,3} - 0.02\, x_1 - 0.67\, x_1^{\,3} + u \\ y = x_1 \end{array}\right\} \quad (9)$$

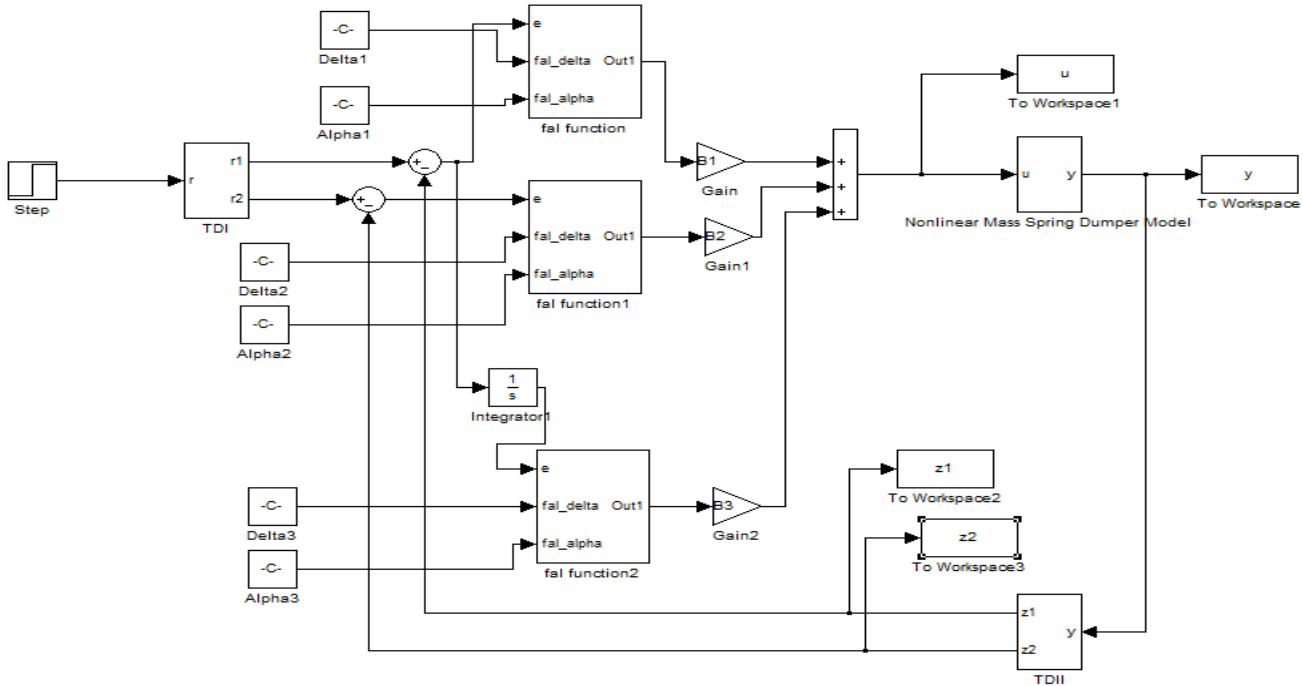

Fig. 3. The Simulink® model for the NPID and the NMSD plant

## V. NUMERICAL SIMULATIONS

The NPID controller based on either the modified Han TD or the proposed INTD and the NMSD mathematical models are designed and numerically simulated using Matlab® /Simulink® as shown in Fig. 3. The values of the parameters for these subsystems are listed in Tables I-III.

TABLE I. THE PARAMETERS OF THE CONTROL LAW

| Parameter | Value |
|---|---|
| $\delta_1$ | 0.1038 |
| $\alpha_1$ | 0.7128 |
| $\beta_1$ | 1.9151 |
| $\delta_2$ | 0.0354 |
| $\alpha_2$ | 0.8680 |
| $\beta_2$ | 2.0130 |
| $\delta_3$ | 1.1916 |
| $\alpha_3$ | 0.9888 |
| $\beta_3$ | 0.0800 |

TABLE II. THE PARAMETERS OF THE MODIFIED HAN TD

| Parameter | Value |
|---|---|
| $R$ | 11.6000 |
| $\delta$ | 0.0005 |

TABLE III. THE PARAMETERS OF THE INTD

| Parameter | Value |
|---|---|
| $\alpha$ | 0.9790 |
| $\beta$ | 5.5872 |
| $\gamma$ | 8.3864 |
| $R$ | 26.5005 |

The numerical simulations are done by using Matlab® ODE45 solver for the models with continuous states. This Runge-Kutta ODE45 solver is a fifth-order method that performs a fourth-order estimate of the error. The reference input to the system is constant linear displacement equals to 0.1 m applied at $t = 0$ sec. The NPID controller is tested for two cases. The numerical simulation of the first testing case is done without adding a measurement noise at the output of the NMSD plant, and the results of this case are shown in Fig.4 and Fig. 5. Also, the numerical results are listed in Table IV.

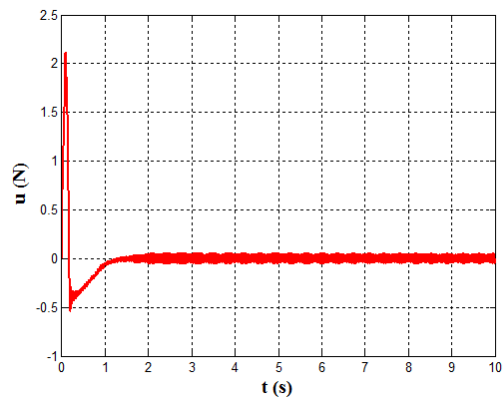

(a)

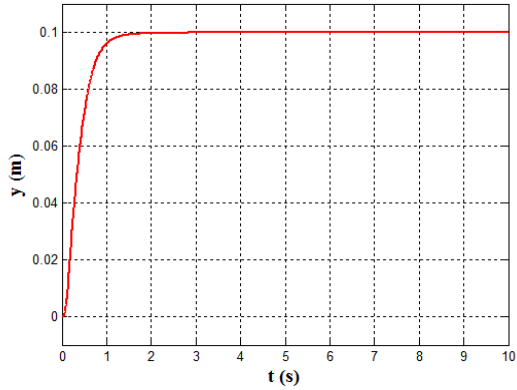

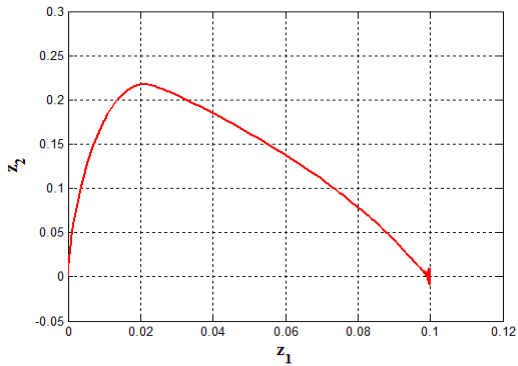

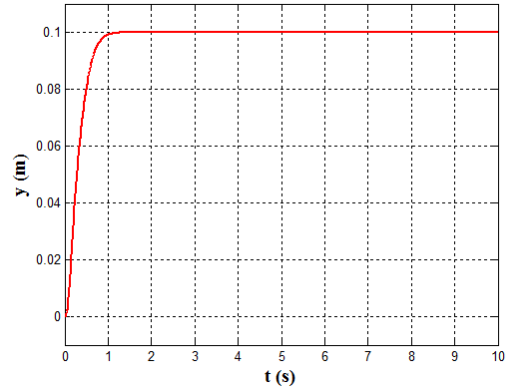

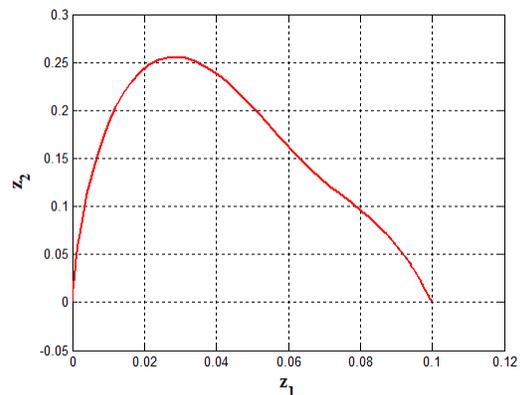

Fig. 4. The simulation results of the NPID based on modified Han TD, (a) The control signal u, (b) The plant output y, and (c) The TD(II) state Trajectory ($z_1, z_2$)

Fig. 5. The simulation results of the NPID based on proposed INTD, (a) The control signal u (b) The plant output y (c) The TD(II) state Trajectory (z1,z2)

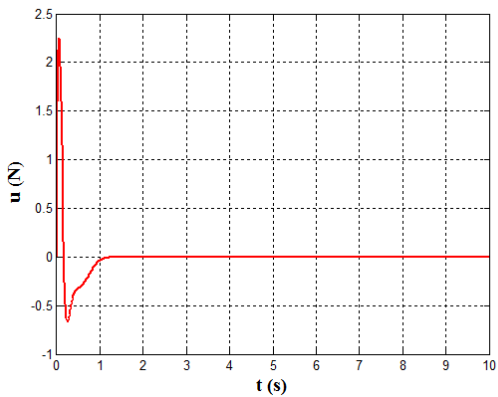

TABLE IV. THE NUMERICAL SIMULATION RESULTS OF CASE 1

| Performance Index | Modified Han TD | Proposed INTD |
|---|---|---|
| IAE | 0.062009 | 0.037965 |
| ITAE | 0.017028 | 0.007961 |
| ITSE | 0.000623 | 0.000325 |
| ISU | 0.532537 | 0.559512 |
| IAU | 0.851211 | 0.540125 |

Where,

$ITAE = \int_0^{10} t \times |r - y| \, dt$ is the integration of the time absolute error for the output signal

$ITSE = \int_0^{10} t \times (r - y)^2 \, dt$ is the integration of the time squared error for the output signal

$IAE = \int_0^{10} |r - y| \, dt$ is the integration of the absolute error for the output signal

$IAU = \int_0^{10} |u| \, dt$ is the integration of absolute of the NPID control signal

$ISU = \int_0^{10} u^2 dt$ is the integration of square of the NPID control signal

Figure 4 (a) shows the chattering in the control signal due to the nonlinear signal of the Han TD [12]. By using the proposed INTD, the chattering in the control signal is significantly reduced (figure 5(a)). The IAU performance index reflects this improvement. The peaking phenomenon previously explained in Lemma (4) appears in the ISU performance index and has the benefit of speeding up the time response of plant output.

The second testing case demonstrated in this work considers adding a measurement noise at the output of the plant. The measurement noise modeled as uniform in the range [-0.001, 0.001] at sampling time 0.001 s. The result of this case shown in Fig. 6, Fig. 7, and Table V.

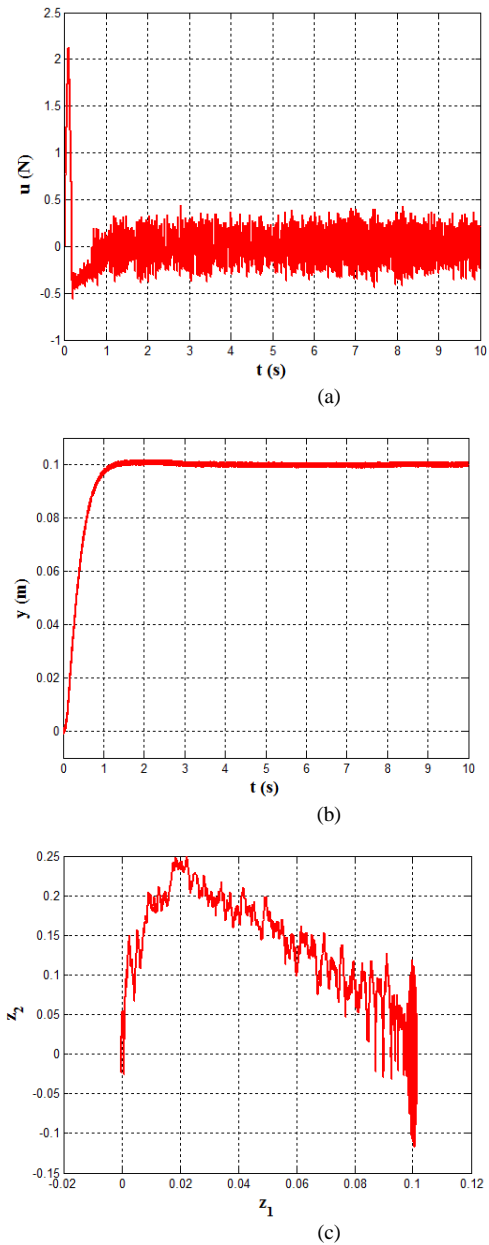

Fig. 6. The simulation results of the NPID based on modified Han TD with measurements noise, (a) The control signal u (b) The plant output y (c) The TD(II) state Trajectory (z1,z2)

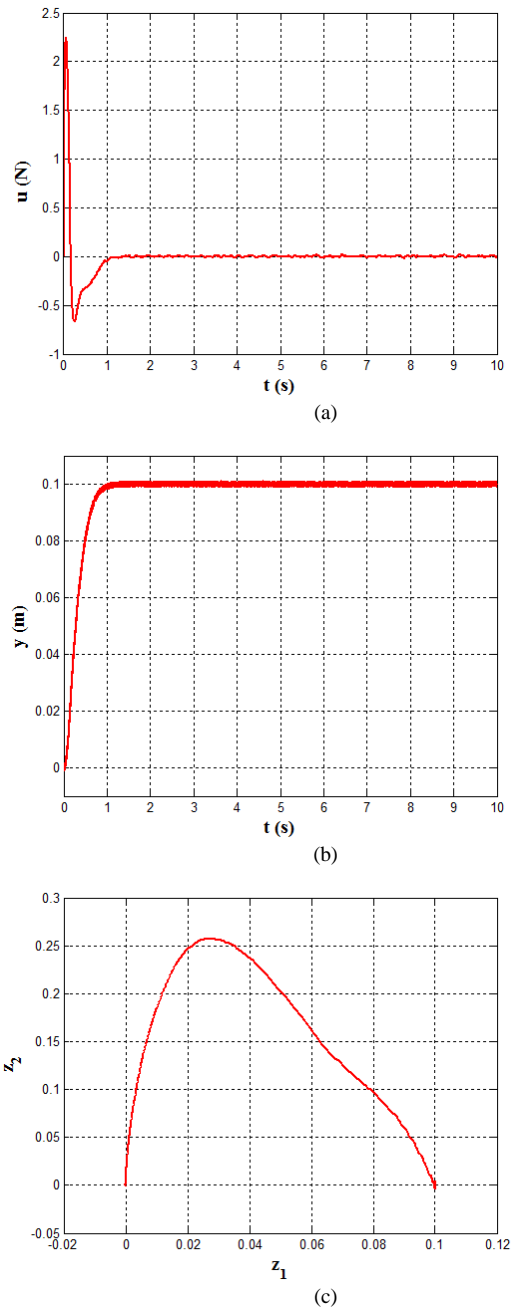

Fig. 7. The simulation results of the NPID based on proposed INTD with measurement noise, (a)The control signal u (b) The plant output y (c) The TD(II) state Trajectory (z1,z2)

TABLE V. THE NUMERICAL SIMULATION RESULTS OF CASE 2

| Performance Index | Modified Han TD | INTD |
|---|---|---|
| IAE | 0.489928 | 0.384329 |
| ITAE | 0.512582 | 0.327142 |
| ITSE | 0.005007 | 0.003458 |
| ISU | 8.449415 | 5.536271 |
| IAU | 26.187044 | 5.561920 |

The band-limiting effect is very clear for the INTD as shown in Fig. 7-(a) with very little fluctuations in the steady state. While the control signal *u* for the case of the TD of [13] is highly affected by the measurement noise (Fig. 6-(a)). Same results are reflected on the output sign*al y* and state-trajectories. Also, the simulations prove that the proposed INTD outperforms the TD offered by [13] with five performance measures as indicated by Table V.

## VI. CONCLUSION

In this article, an improved type of nonlinear tracking differentiator is developed to obtain higher derivatives of reference signal to achieve tracking with high robustness against measurement noise. The proposed tracking differentiator is proven to be globally asymptotically stable. It converges to the exact derivatives of the signal independent of the initial differentiation error. The INTD has an under damped effect which lead directly to peaking phenomenon. Knowing that the INTD is a continuous structure which comprises of rectilinear and non-linear parts, the noise and chattering phenomenon has been reduced adequately, the reason is due ot the high fidelity that the INTD has when generting the derivatives of the signal. Also dynamical performance are enhanced apparently. The Simulation experiments show the feasibility of integrating the proposed INTD with the nonlinear combinations of the error profile to design a nonlinear PID controller for MSD system which can be considered as an alternative and efficient control method to solve real control design for such nonlinear systems. The new configuration with the proposed INTD achieves fast arrival and smooth tracking to the input signal. Finally, the performance of the nonlinear MSD system has been enhanced dramatically.

As a future work, the value of *R* in this work may be varied and based on certain adaptive law, the optimal value for *R* may be chosen which makes the TD producing better results.

ACKNOWLEDGMENT

The authors appreciate the electrical engineering department laboratories for direct help and support to finish this research.